\documentclass[preprint,12pt]{elsarticle}
\usepackage{float}
\usepackage[utf8]{inputenc}
\usepackage[T1]{fontenc}
\usepackage{amsmath,amssymb,amsthm,mathtools,bm}
\usepackage{enumitem,booktabs,array}
\usepackage{tikz}
\usetikzlibrary{arrows.meta}
\usepackage[expansion=false]{microtype}
\usepackage{placeins}
\usepackage[hidelinks]{hyperref}
\usepackage{xurl}
\makeatletter
\def\ps@pprintTitle{\let\@oddhead\@empty\let\@evenhead\@empty\let\@oddfoot\@empty\let\@evenfoot\@empty}
\makeatother
\hypersetup{pdftitle={Fuzzy vectors from first and second moments: Affine consistency, carrying and sparsest realizations},pdfauthor={Jan Schneider; Paul Schneider}}
\theoremstyle{plain}
\newtheorem{theorem}{Theorem}[section]
\newtheorem{proposition}[theorem]{Proposition}
\newtheorem{lemma}[theorem]{Lemma}

\theoremstyle{definition}
\newtheorem{definition}[theorem]{Definition}
\newtheorem{example}[theorem]{Example}
\newtheorem{remark}[theorem]{Remark}
\newcommand{\R}{\mathbb R}
\newcommand{\E}{\mathbb E}
\newcommand{\Prob}{\mathbb P}
\newcommand{\F}{\mathbb F}

\newcommand{\Sspace}{\mathcal S}

\newcommand{\Sphere}{\mathbb S}
\newcommand{\Ell}{\operatorname{Ell}}
\newcommand{\Img}{\operatorname{Im}}
\newcommand{\rank}{\operatorname{rank}}
\newcommand{\aff}{\operatorname{aff}}
\newcommand{\conv}{\operatorname{conv}}
\newcommand{\supp}{\operatorname{supp}}
\newcommand{\Cov}{\operatorname{Cov}}
\newcommand{\Var}{\operatorname{Var}}
\newcommand{\tr}{\operatorname{tr}}
\newcommand{\law}{\mathcal L}
\newcommand{\GL}{\operatorname{GL}}
\newcommand{\Orth}{\operatorname O}
\newcommand{\Unif}{\operatorname{Unif}}

\newcommand{\bA}{\bm{A}}
\newcommand{\bB}{\bm{B}}
\newcommand{\bC}{\bm{C}}
\newcommand{\bG}{\bm{G}}
\newcommand{\bI}{\bm{I}}
\newcommand{\bJ}{\bm{J}}
\newcommand{\bL}{\bm{L}}

\newcommand{\bPi}{\bm{\Pi}}
\newcommand{\bQ}{\bm{Q}}
\newcommand{\bS}{\bm{S}}
\newcommand{\bSigma}{\bm{\Sigma}}
\newcommand{\bTheta}{\bm{\Theta}}
\newcommand{\bX}{\bm{X}}
\newcommand{\bZ}{\bm{Z}}
\newcommand{\ba}{\bm{a}}
\newcommand{\bb}{\bm{b}}
\newcommand{\bc}{\bm{c}}
\newcommand{\be}{\bm{e}}
\newcommand{\bmu}{\bm{\mu}}
\newcommand{\bu}{\bm{u}}
\newcommand{\bv}{\bm{v}}
\newcommand{\bx}{\bm{x}}
\newcommand{\by}{\bm{y}}
\newcommand{\bz}{\bm{z}}
\newcommand{\bzeta}{\bm{\zeta}}
\newcommand{\norm}[1]{\left\lVert#1\right\rVert}
\newcommand{\dd}{\,\mathrm d}
\providecommand{\doi}[1]{\href{https://doi.org/#1}{doi:\nolinkurl{#1}}}

\setlist{leftmargin=*}
\allowdisplaybreaks[1]

\begin{document}

\begin{frontmatter}
\title{Fuzzy vectors from first and second moments:\\ Affine consistency, carrying and sparsest realizations}
\author[wro]{Jan Schneider}
\author[usi,sfi]{Paul Schneider\corref{cor1}}
\ead{paul.schneider@usi.ch}
\cortext[cor1]{Corresponding author.}
\affiliation[wro]{organization={Wroc{\l}aw University of Science and Technology},city={Wroc{\l}aw},country={Poland}}
\affiliation[usi]{organization={Universit\`a della Svizzera italiana (USI)},city={Lugano},country={Switzerland}}
\affiliation[sfi]{organization={Swiss Finance Institute (SFI)},country={Switzerland}}

\begin{abstract}\small

When the probability law of a random vector cannot be identified, its mean and possibly singular covariance can often still be elicited. We consider fuzzy vectors with compact convex cuts as mathematically tractable representations of such a moment pair, and ask how they must be designed to rely only on this information.
We postulate two modeling desiderata, independence of the choice of origin and basis, and invariance under embeddings into a higher-dimensional space. Under these requirements alone we show that, in every ambient dimension, every level set of the fuzzy vector must be a Mahalanobis ellipsoid centered at the mean.
Its complete specification then amounts to the choice of a radial profile for each covariance rank. We call a fuzzy vector a carrier of the moment pair if at least one of its random-set selections has exactly the given mean and covariance. For an ellipsoidal fuzzy vector of rank \(p\ge1\) with radial profile \(\rho\), carrying is equivalent to \(\int_0^1\rho(\alpha)^2\dd\alpha\ge p\), and equality characterizes the inclusion-minimal carriers within the fixed-shape Mahalanobis class.
A discrete random vector with the same mean and covariance and the fewest equally weighted atoms has its atoms on a regular simplex in Mahalanobis coordinates. Completing this simplex under the rotations left free by the moments yields the minimum-volume crisp ellipsoid in the state space that carries them. The ellipsoidal representation also encodes interaction through the covariance matrix.
\end{abstract}
\end{frontmatter}

\section{Introduction}\label{sec:intro}
Uncertain multivariate quantities arise throughout physics, finance, engineering, and the social sciences. In many such settings the underlying probability law cannot be identified reliably, whereas its mean vector and covariance matrix can be estimated from data, assessed by experts, or reported as summary statistics. The covariance may be singular, for instance when there are more variables than observations. These two moments describe location, marginal dispersion and linear dependence, but they do not determine a probability law, since infinitely many distributions may share the same moment pair. We ask what fuzzy-vector representation can be constructed from this information without introducing distinctions unwarranted by it.

Our starting point is a principle of \emph{informational restraint}. A representation based only on specified information should not depend on arbitrary features of its description. We formulate this principle as \emph{coordinate neutrality} and \emph{ambient neutrality}. The former requires a representation to transform with the moments under invertible affine changes of coordinates. The latter requires that adding deterministic coordinates leave the intrinsic uncertainty unchanged. Their joint mathematical expression we term \emph{affine consistency}.

The direction of the argument is important. We do not first choose ellipsoids and then verify their invariance. We specify the admissible fuzzy vectors and the consistency requirements, and characterize every construction that satisfies them. Once these premises are accepted, the geometry is a consequence rather than a design choice. Theorem~\ref{thm:characterization} shows that every cut must be a Mahalanobis ellipsoid in the state space, regardless of ambient dimension or singularity. Its radius may depend on the membership level and the covariance rank, but not on the choice of origin and basis. We term the function assigning to each membership level the Mahalanobis radius of the corresponding cut \emph{radial profile}. The remaining freedom is therefore exactly one radial profile for each positive rank. Section~\ref{sec:role-assumptions} explains the separate roles of the assumptions.

We then introduce an additional compatibility requirement, not implied by neutrality. A fuzzy vector \emph{carries} the prescribed moments if one of its classical random-set selections has exactly the given mean and covariance. The notion of selection is standard, whereas carrying is the new concept introduced here for the moment-representation problem.

Theorem~\ref{thm:carrying} gives an exact integral condition on the radial profile of the fuzzy vector, and equality in it characterizes the minimal carriers within the fixed-shape Mahalanobis class. A bivariate example shows that minimal carriers are not unique, since different profiles satisfy the same condition. Combined with the behavior of ellipsoids under projections, the carrying condition also shows that extending affine consistency to rank-reducing maps is incompatible with carrying moment pairs of every rank, and with normalization, that is, equality in this condition, at two different ranks (Proposition~\ref{prop:rankreduction-carrying}). The restriction to injective maps is therefore substantive.

There is a complementary discrete perspective. Any covariance of rank \(p\) can be realized exactly by a discrete probability law on \(p+1\) equally weighted points, and \(p+1\) is the minimum, because the points must span the \(p\)-dimensional state space. Such points can be computed efficiently \cite{FastEmpirical}. Theorem~\ref{thm:simplex} shows that, after whitening, every set of \(p+1\) equally weighted points with the prescribed mean and covariance is a regular simplex. Its orientation is arbitrary, since every rotation of a realization is again a realization. Taking the union of all rotated copies yields precisely the crisp Mahalanobis ellipsoid of radius \(\sqrt p\), the unique minimum-volume crisp ellipsoid in the state space that carries the moments. The continuous ellipsoidal geometry and the sparsest discrete geometry thus meet through symmetry completion. Symmetry  singles out the ellipsoid, since fuzzy vectors whose cuts are dilated copies of such a simplex carry the moments with cuts that are proper subsets of the ellipsoidal ones.

Section~\ref{sec:interaction} shows how the same ellipsoidal representation encodes interaction. The spread of a linear combination of the coordinates is governed by the covariance quadratic form. It is never larger than under the minimum-based non-interactive joint of the same marginals, that is, of the projections of the joint onto the coordinates, and it can be strictly smaller even for uncorrelated coordinates.

\paragraph{Our contributions} We offer the following: the characterization of all moment-based assignments, including singular covariances and consistency across ambient dimensions, which leaves exactly one radial profile for each intrinsic rank, the notion of carrying with its exact criterion and its minimal carriers of fixed shape, the obstruction to consistency under rank-reducing maps, and the identification of the radius-\(\sqrt p\) Mahalanobis ellipsoid as the unique minimum-volume carrier in the state space, which the rotational completion of a sparsest realization recovers.

\paragraph{Related work}
Turning partial probabilistic knowledge into possibility distributions is a long-standing theme of possibility theory. Dubois, Foulloy, Mauris and Prade \cite{DFMP04} develop probability--possibility transformations based on probabilistic inequalities, including the case of a prescribed mean and variance. Baudrit and Dubois \cite{BaudritDubois2006} compare representations of incomplete probabilistic knowledge. Universal transformations require domination of every compatible law. Carrying is an existential counterpart, built on the standard notion of random-set selection (Section~\ref{sec:selection}). Related foundations include selection-based means and variances of imprecise quantities \cite{DuboisPrade1987,Kruse1987}, random sets as ill-perceived random variables \cite{CousoDuboisSanchez2014}, random-set theory \cite{Molchanov2017,Artstein1983}, and identification through selections in econometrics \cite{BeresteanuMolchanovMolinari2011}. Couso and Dubois \cite{CousoDubois2009} contrast a crisp variance of a fuzzy random variable with a fuzzy interval of variances of the classical random variables it may represent, and propose an interval-valued variance within imprecise probabilities. Salamanca and Couso \cite{SalamancaCouso2022} give necessary and sufficient conditions for a random closed set in Euclidean space to admit a selection of minimal Kruse variance. Carrying asks instead for a selection whose covariance equals a prescribed matrix. A different line defines possibilistic means and variances directly from the cuts of a fuzzy number \cite{CarlssonFuller2001,FullerMajlender2003}.

Ellipsoidal fuzzy objects and arguments favoring them predate the present construction. Kreinovich, Beck and Nguyen \cite{KreinovichBeckNguyen2007} motivate ellipsoid-shaped fuzzy sets through quadratic approximation and elicitation, and develop their use in data processing. Neumaier's potential clouds \cite{Neumaier2004} separate a chosen potential from scalar level maps, and quadratic potentials include covariance-shaped ellipsoids. Wang, Shi, Wang and Zhang \cite{WangShiWangZhang2014} define fuzzy \(n\)-ellipsoid numbers by prescribing ellipsoidal cuts and establish cut-representation results for that class. Earlier, Celmi\c{n}\v{s} \cite{Celmins1987} fits models to fuzzy vector data with conical membership functions, whose cuts are ellipsoids determined by a positive definite matrix. Finkelstein, Kosheleva and Kreinovich \cite{FinkelsteinKoshelevaKreinovich1996} characterize ellipsoids among minimal-parameter convex families selected by affine-invariant optimality criteria with a unique optimum. In multivariate statistics, Zuo and Serfling \cite{ZuoSerfling2000}, generalizing Liu and Singh \cite{LiuSingh1993}, show that for an elliptically distributed random vector every affine-invariant depth function attaining its maximum at the center has Mahalanobis-ellipsoidal depth regions. The underlying full-rank rotational-symmetry argument is therefore well established, and we do not claim it as new. Hu\v{s}ek \cite{Husek2012} describes uncertain parameters of linear dynamical systems by multidimensional fuzzy sets with ellipsoidal cuts. Closest to the present setting, Lesniewska-Choquet, Mauris, Atto and Mercier \cite{LesniewskaChoquet2020} define elliptical possibility distributions by a probability--possibility transformation of elliptical probability laws. Here, by contrast, no elliptical probability law is assumed. Our question concerns the assignment rule itself. Which mappings from moment pairs to fuzzy vectors satisfy the two neutrality requirements when no ellipsoidal geometry is assumed in advance? Our characterization includes singular covariances and identifies the dependence on intrinsic rank across ambient dimensions.

Support-function representations of fuzzy vectors go back to Puri and Ralescu \cite{PuriRalescu1985,PuriRalescu1986} and have since been developed in several forms, see for example \cite{LiangNavaraVetterlein2009}. We use the levelwise support-function representation both to describe affine images and in the proof of the characterization theorem.

Arithmetic of interactive fuzzy numbers through a joint possibility distribution goes back to Dubois and Prade \cite{DuboisPrade1981}. Interactive fuzzy numbers and joint possibility distributions are studied by Full\'er and Majlender \cite{FullerMajlender2004} and Carlsson, Full\'er and Majlender \cite{CarlssonFullerMajlender2005}. Related multivariate models include oblique fuzzy vectors \cite{InuiguchiRamikTanino2003} and quadratic or exponential possibility distributions with covariance-type shape matrices \cite{InuiguchiSakawa1995,TanakaIshibuchiYoshikawa1995,TanakaGuo1998,HoussineauBishop2018}.

Point sets of smallest size with prescribed moments occur in the truncated moment problem \cite{CurtoFialkow2005,Lasserre2010} and as sigma points in nonlinear filtering \cite{Julier2003,JulierUhlmann2004}. Equally weighted integration formulas of degree two with \(n+1\) nodes go back to Stroud \cite{Stroud1957,Stroud1960}, and Xiu \cite{Xiu2008} shows that, after affine normalization, such nodes form a regular simplex. Tyler \cite{Tyler2010} shows that for at most \(p+1\) points in general position in \(\R^p\) the only affine-equivariant location statistic is the sample mean, and every affine-equivariant scatter matrix is a fixed multiple of the sample covariance. Part~(i) of Theorem~\ref{thm:simplex} thus has direct antecedents, and its short proof is included for completeness. Its role here is to connect sparsest realizations with carrying and rotational completion. The covariance scenarios of Multerer, Schneider and Sen \cite{FastEmpirical} compute equally weighted moment realizations with a single Householder reflection. In this work their role is to reveal the relation between sparsest discrete realizations and the symmetry imposed by the axioms.

\smallskip

\paragraph{Paper organization} The paper is organized as follows. Section~\ref{sec:preliminaries} introduces the three equivalent representations of fuzzy vectors, affine images, moment geometry, selections, and carrying. Section~\ref{sec:affine} characterizes affine-consistent constructions, describes their affine images and linear combinations, and examines the role of the assumptions. Section~\ref{sec:carrying} derives the carrying condition on the radial profile, Section~\ref{sec:simplex} relates it to sparsest discrete realizations and symmetry completion, and Section~\ref{sec:discussion} discusses limitations and further work.

\section{Preliminaries and basic concepts}\label{sec:preliminaries}

\subsection{Fuzzy vectors}

As it is convenient for the subsequent analysis, we shall define a fuzzy vector in three equivalent ways, by its membership function, by its nested family of level sets, and by the support functions of those level sets. We start with the membership function.

\begin{definition}[Fuzzy vector]\label{def:fuzzy}
A \(d\)-dimensional \emph{fuzzy vector} \(\xi\) is specified by a normal, upper semicontinuous membership function
\[
   \pi_\xi:\R^d\longrightarrow[0,1]
\]
with bounded support and convex positive level sets. Equivalently, for every \(\alpha\in(0,1]\), the \(\alpha\)-cut
\[
   C_\alpha(\xi)
      :=
   \{\bx\in\R^d:\pi_\xi(\bx)\ge\alpha\}
\]
is nonempty, compact and convex, and all cuts are contained in a common bounded set, the fuzzy vector's support
\[
   C_0(\xi)
      :=
   \overline{\{\bx\in\R^d:\pi_\xi(\bx)>0\}},
\]
the \(0\)-cut.

We denote the class of such fuzzy vectors by \(\F_d\). We call a fuzzy vector \emph{crisp} if its membership function is the indicator function of a single nonempty compact convex set.
\end{definition}

Here, normal means that the membership function attains the value \(1\). Upper semicontinuity is equivalent to closedness of all positive cuts. The adjective \emph{crisp} therefore includes, but is not restricted to, singleton-valued fuzzy vectors.

\smallskip
\paragraph{Level-set representation}
The cuts of a fuzzy vector are nested and satisfy
\begin{equation}\label{eq:cutregularity}
 C_\alpha(\xi)
   =
 \bigcap_{0<\gamma<\alpha}C_\gamma(\xi),
 \qquad 0<\alpha\le1,
 \qquad
 C_0(\xi)
   =
 \overline{\bigcup_{\alpha>0}C_\alpha(\xi)}.
\end{equation}
Conversely, every uniformly bounded nested family of nonempty compact convex sets satisfying \eqref{eq:cutregularity} determines a unique fuzzy vector through
\[
   \pi_\xi(\bx)
      =
   \sup\{\alpha\in[0,1]:\bx\in C_\alpha(\xi)\},
   \qquad \sup\emptyset=0.
\]
Thus the membership-function and level-set descriptions are equivalent.

\smallskip
\paragraph{Support-function representation}
For a nonempty compact convex set \(C\subset\R^d\), its support function is
\[
   h_C(\bu)
      :=
   \max_{\bx\in C}\bu^\top \bx,
   \qquad \bu\in\R^d.
\]
The set \(C\) is uniquely characterized by \(h_C\), since
\[
   C
      =
   \bigcap_{\bu\in\mathbb S^{d-1}}
   \{\bx\in\R^d:\bu^\top \bx\le h_C(\bu)\}.
\]
Accordingly, a fuzzy vector \(\xi\) is equivalently represented by its levelwise support-function map
\begin{equation}\label{eq:fuzzysupport}
   s_\xi(\alpha,\bu)
      :=
   h_{C_\alpha(\xi)}(\bu),
   \qquad
   (\alpha,\bu)\in[0,1]\times\mathbb S^{d-1}.
\end{equation}
Indeed, \eqref{eq:fuzzysupport} recovers every cut, and hence the membership function. Characterizations of the functions \(s\) that arise in this way are standard in the support-function theory of fuzzy vectors \cite{LiangNavaraVetterlein2009}. This representation is used extensively in \cite{SU1,SU2}. It also turns Minkowski addition of cuts into pointwise addition of their support functions and identifies the Hausdorff metric on compact convex cuts with the corresponding sup norm.

\subsection{Affine transformations}

Let
\[
   T(\bx)=\bA\bx+\bb:\R^d\longrightarrow\R^m.
\]
The image of a fuzzy vector \(\xi\) transformed by \(T\), denoted \(T\cdot\xi\), is defined by Zadeh's extension principle,
\[
   \pi_{T\cdot\xi}(\by)
      =
   \sup\{\pi_\xi(\bx):T(\bx)=\by\},
   \qquad \sup\emptyset=0.
\]
For the class of fuzzy vectors considered here this is equivalent to the cutwise identity
\begin{equation}\label{eq:cutimage}
   C_\alpha(T\cdot\xi)
      =
   T C_\alpha(\xi),
   \qquad 0\le\alpha\le1,
\end{equation}
including the support cut, compare \cite{Nguyen1978}.

\smallskip
\paragraph{Support functions under affine maps}
If \(C\subset\R^d\) is nonempty, compact and convex, then
\begin{equation}\label{eq:supportaffine}
   h_{TC}(\bu)
      =
   h_C(\bA^\top \bu)+\bu^\top \bb,
   \qquad \bu\in\R^m.
\end{equation}
Thus affine transformations of fuzzy vectors become simple transformations of their levelwise support functions.

\subsection{Moment geometry}

Let \(\bmu\in\R^d\) and let \(\bSigma\in\R^{d\times d}\) be symmetric positive semidefinite. Write
\[
   p:=\rank\bSigma,
   \qquad
   \Sspace(\bmu,\bSigma)
      :=
   \bmu+\Img\bSigma.
\]
We call \(p\) the \emph{intrinsic rank} of the moment pair and \(\Sspace(\bmu,\bSigma)\) its \emph{state space}. If \(\bSigma\) is nonsingular, then \(\Sspace(\bmu,\bSigma)=\R^d\). If \(\bSigma\) is singular, the state space is a proper affine subspace of \(\R^d\).

\smallskip
\paragraph{Mahalanobis radius}
Let \(\bSigma^+\) denote the Moore--Penrose inverse of \(\bSigma\), and denote
\[
   \bSigma^{+/2}:=(\bSigma^+)^{1/2}.
\]
The \emph{Mahalanobis radius} associated with \((\bmu,\bSigma)\) is
\begin{equation}\label{eq:radius}
   t_{\bmu,\bSigma}(\bx)
      :=
   \bigl\|\bSigma^{+/2}(\bx-\bmu)\bigr\|,
   \qquad
   \bx\in\Sspace(\bmu,\bSigma),
\end{equation}
and we set \(t_{\bmu,\bSigma}(\bx)=+\infty\) outside the state space. When the moment pair is implicit from the context, we write simply \(t(\bx)\).

For \(\bx\in\Sspace(\bmu,\bSigma)\) we call \(\bz:=\bSigma^{+/2}(\bx-\bmu)\in\Img\bSigma\) the \emph{Mahalanobis coordinates} of \(\bx\). The map \(\bx\mapsto\bz\) is a bijection from \(\Sspace(\bmu,\bSigma)\) onto \(\Img\bSigma\) with inverse \(\bz\mapsto\bmu+\bSigma^{1/2}\bz\), and \(t(\bx)=\|\bz\|\). A random vector with moments \((\bmu,\bSigma)\) lies in \(\Sspace(\bmu,\bSigma)\) almost surely, as shown below, and has Mahalanobis coordinates with mean \(0\) and covariance \(\bPi_{\bSigma}\), where \(\bPi_{\bSigma}:=\bSigma\bSigma^+\) is the orthogonal projector onto \(\Img\bSigma\).

\smallskip
\paragraph{Mahalanobis ellipsoids}
Let
\[
   \mathbb B_{\Img\bSigma}
      :=
   \{\bz\in\Img\bSigma:\|\bz\|\le1\}.
\]
Its boundary is the intrinsic unit sphere \(\Sphere_{\Img\bSigma}:=\{\bz\in\Img\bSigma:\|\bz\|=1\}\), and, for \(p\ge1\), \(\Unif(\Sphere_{\Img\bSigma})\) denotes the normalized surface measure on it. For \(\Img\bSigma=\R^p\) we write \(\mathbb B_p\) for the closed unit ball.
For \(r\ge0\), the sublevel set
\begin{equation}\label{eq:mahalanobisball}
   E_r(\bmu,\bSigma)
      :=
   \{\bx:t_{\bmu,\bSigma}(\bx)\le r\}
      =
   \bmu+r\bSigma^{1/2}\mathbb B_{\Img\bSigma}
\end{equation}
is the \emph{Mahalanobis ellipsoid of radius \(r\)} associated with the moment pair. When \(p<d\), this ellipsoid is degenerate as a subset of \(\R^d\), but for \(r>0\) and \(p\ge1\) it is a full \(p\)-dimensional ellipsoid in the state space.

\smallskip
\paragraph{Support implied by the moments}
If a random vector \(\bX\) satisfies
\[
   \E \bX=\bmu,
   \qquad
   \Cov \bX=\bSigma,
\]
then \(\bX\in\Sspace(\bmu,\bSigma)\) almost surely and
\begin{equation}\label{eq:affinesupport}
   \aff\supp\law(\bX)
      =
   \Sspace(\bmu,\bSigma).
\end{equation}
To see the first statement, if \(\bv\in\ker\bSigma\), then
\[
   \Var(\bv^\top \bX)=\bv^\top\bSigma \bv=0,
\]
so \(\bv^\top \bX=\bv^\top\bmu\) almost surely. If the affine hull of the support were a proper subspace of \(\Sspace(\bmu,\bSigma)\), there would be an additional nonzero zero-variance direction. This elementary observation is classical \cite{Rao1973}.

\subsection{Selections}\label{sec:selection}

A fuzzy vector is connected with ordinary probability by drawing a membership level at random and then a point in the corresponding cut. Throughout, \(U\sim\Unif(0,1)\) denotes the random level. Since the cuts of \(\xi\) are nonempty, compact and nested and satisfy \eqref{eq:cutregularity}, the event \(\{C_U(\xi)\cap K\neq\emptyset\}\) equals \(\{U\le\max_{\bx\in K}\pi_\xi(\bx)\}\) for every nonempty compact \(K\), the maximum being attained by upper semicontinuity of \(\pi_\xi\). Hence \(C_U(\xi)\) is a random compact set \cite{Molchanov2017}.

\begin{definition}[Selection]\label{def:selection}
Let \(U\sim\Unif(0,1)\), and let \(\bX\) be a random vector in \(\R^d\) defined on the same probability space. Then \(\bX\) is a \emph{selection} of \(\xi\) if
\[
   \bX\in C_U(\xi)
   \qquad\text{almost surely}.
\]
\end{definition}
In random-set terminology \(\bX\) is a selection of the random set \(C_U(\xi)\) \cite{Aumann1965,Artstein1983,Molchanov2017}. Only the joint law of \((U,\bX)\) matters. In particular, \(\bX\) may use randomness beyond \(U\) to choose a point within the realized cut.

\smallskip
\paragraph{Levelwise kernels}
Selections can be described equivalently by their conditional laws given the level. Let \(\alpha\mapsto P_\alpha\) be a measurable probability kernel from \((0,1)\) to \(\R^d\) such that
\begin{equation}\label{eq:kernel-support}
   P_\alpha\bigl(C_\alpha(\xi)\bigr)=1
   \qquad\text{for almost every }\alpha\in(0,1).
\end{equation}
If \((U,\bX)\) has joint law \(\dd\alpha\otimes P_\alpha(\dd \bx)\) on \((0,1)\times\R^d\), then \(\bX\) is a selection of \(\xi\), and its law is the uniform mixture
\begin{equation}\label{eq:kernel-mixture}
   P(A)
      :=
   \int_0^1 P_\alpha(A)\,\dd\alpha,
   \qquad A\in\mathcal B(\R^d).
\end{equation}
Conversely, every selection \(\bX\) has a regular conditional law \(P_\alpha\) given \(U=\alpha\), because \(\R^d\) is Polish. Since \(\bX\in C_U(\xi)\) almost surely, this kernel satisfies \eqref{eq:kernel-support}. The laws of selections are therefore exactly the mixtures \eqref{eq:kernel-mixture} of kernels satisfying \eqref{eq:kernel-support}. The kernels need not have Lebesgue densities. They may be singular with respect to Lebesgue measure on \(\aff C_\alpha(\xi)\), for instance point masses or laws concentrated on the boundary of a cut.

\smallskip
\paragraph{Coverage}
The random set \(C_U(\xi)\) also recovers the membership function. For every fixed \(\bx\in\R^d\),
\[
   \bx\in C_U(\xi)
   \quad\Longleftrightarrow\quad
   U\le\pi_\xi(\bx),
\]
and therefore
\begin{equation}\label{eq:coverage}
   \Prob\bigl(\bx\in C_U(\xi)\bigr)
      =
   \Prob\bigl(U\le\pi_\xi(\bx)\bigr)
      =
   \pi_\xi(\bx).
\end{equation}
Thus \(\pi_\xi\) is the \emph{one-point coverage function} of the random set \(C_U(\xi)\). This consonant random-set representation of fuzzy sets is classical \cite{Goodman1982,Molchanov2017}.

\subsection{Carrying prescribed moments}
We now introduce a new and separate compatibility concept for the present moment-based problem.

\begin{definition}[Carrying prescribed moments]\label{def:carrying}
A fuzzy vector \(\xi\) \emph{carries} the moment pair \((\bmu,\bSigma)\) if it admits a selection \(\bX\) such that
\[
   \E \bX=\bmu,
   \qquad
   \Cov \bX=\bSigma.
\]
\end{definition}

Equivalently, \(\xi\) carries \((\bmu,\bSigma)\) if there is a kernel satisfying \eqref{eq:kernel-support} for which the mixture \eqref{eq:kernel-mixture} has these first two moments. Carrying is therefore an existential compatibility requirement. At least one probability law compatible with the fuzzy vector must realize the given moment pair. By \eqref{eq:affinesupport}, every carrier necessarily satisfies
\begin{equation}\label{eq:carrierstate}
   \Sspace(\bmu,\bSigma)
      \subseteq
   \aff C_0(\xi).
\end{equation}

\smallskip
\paragraph{Existential and universal compatibility}
Carrying asks for one compatible law. By Artstein's theorem \cite{Artstein1983,Molchanov2017}, a probability law \(P\) on \(\R^d\) is the law of a selection of \(\xi\) if and only if
\[
   P(K)\le\max_{\bx\in K}\pi_\xi(\bx)
   \qquad\text{for every nonempty compact }K\subset\R^d,
\]
that is, if and only if \(P\) is dominated by the possibility measure of \(\xi\). Hence \(\xi\) carries \((\bmu,\bSigma)\) precisely when its possibility measure dominates at least one law with mean \(\bmu\) and covariance \(\bSigma\). This differs from probability--possibility transformations that seek a possibility distribution dominating every probability law satisfying the available probabilistic constraints \cite{DFMP04,BaudritDubois2006}.

Carrying ensures that the moments are realizable within the random-set interpretation of the chosen fuzzy vector. It does not assert that every law sharing those moments, or the actual data-generating law, is a selection law. For example, a bounded carrier cannot have a Gaussian law with the same nonzero covariance among its selection laws. Its cuts are therefore not distribution-free confidence or prediction regions. This distinction belongs to the modeling interpretation and does not qualify the mathematical results.

\section{Affine consistency}\label{sec:affine}

We now consider constructions that assign a fuzzy vector to every moment pair $\bmu$ and $\bSigma$.
Such a construction is a family of maps
\[
   \Phi_d:\mathcal M_d\longrightarrow\F_d,
   \qquad
   \mathcal M_d
      :=
   \R^d\times\R_{\succeq0}^{d\times d},
   \qquad
   (\bmu,\bSigma)\longmapsto\Phi_d(\bmu,\bSigma).
\]

This is where a type of informational restraint enters (compare with \cite{DuboisPrade1988}).
A construction based only on the moment pair should depend on the information the moments carry, not on how that information is described. In particular, it should not depend on the choice of origin and basis in which \((\bmu,\bSigma)\) is expressed, nor on the number of deterministic coordinates in which the uncertain quantity is embedded. We formulate this principle through two requirements.

\begin{definition}[Affine consistency]\label{def:affine}
The family \(\Phi=(\Phi_d)_{d\ge1}\) of maps \(\Phi_d:\mathcal M_d\to\F_d\) is \emph{affine-consistent} if it satisfies the following two conditions.

\begin{enumerate}[label=(\roman*)]
\item \label{prop:coordinateneut} \emph{Coordinate neutrality.}
For every invertible affine map \(T(\bx)=\bA\bx+\bb\) of \(\R^d\), with \(\bA\in\GL(d)\) and \(\bb\in\R^d\),
\begin{equation}\label{eq:equivariance}
   \Phi_d(\bA\bmu+\bb,\bA\bSigma \bA^\top)
      =
   T\cdot\Phi_d(\bmu,\bSigma).
\end{equation}

\item \label{prop:ambientneut}\emph{Ambient neutrality.}
For every \(m\ge d\), let \(\iota_{d,m}:\R^d\hookrightarrow\R^m\), \(\iota_{d,m}(\bx)=(\bx,0)\), be the canonical embedding. Then
\begin{equation}\label{eq:embedding}
   \Phi_m\!\left(
      \begin{pmatrix}\bmu\\0\end{pmatrix},
      \begin{pmatrix}\bSigma&0\\0&0\end{pmatrix}
   \right)
      =
   \iota_{d,m}\cdot\Phi_d(\bmu,\bSigma).
\end{equation}
\end{enumerate}
\end{definition}

We state affine consistency through the two conditions of Definition~\ref{def:affine} because each makes our two requirements of informational restraint above transparent. Mathematically, they can, however, be combined into a single condition.

\begin{lemma}[Injective affine consistency]\label{lem:injective}
A family \(\Phi=(\Phi_d)_{d\ge1}\) of maps \(\Phi_d:\mathcal M_d\to\F_d\) is affine-consistent if and only if, for every injective affine map \(T:\R^d\to\R^m\), \(T(\bx)=\bA\bx+\bb\),
\begin{equation}\label{eq:injective}
   \Phi_m(\bA\bmu+\bb,\bA\bSigma \bA^\top)
      =
   T\cdot\Phi_d(\bmu,\bSigma).
\end{equation}
\end{lemma}

\begin{proof}
Taking \(m=d\) in \eqref{eq:injective} gives coordinate neutrality, and taking \(T=\iota_{d,m}\) gives ambient neutrality. Conversely, complete the columns of \(\bA\) to a basis of \(\R^m\). This gives \(\bB\in\GL(m)\) with \(\bA=\bB\iota_{d,m}\), so that \(T=S\circ\iota_{d,m}\) with \(S(\by)=\bB\by+\bb\). Images under the extension principle compose, that is, \((S\circ\iota_{d,m})\cdot\xi=S\cdot(\iota_{d,m}\cdot\xi)\) for every fuzzy vector \(\xi\). Applying ambient neutrality and then coordinate neutrality yields \eqref{eq:injective}.
\end{proof}

The restriction to injective maps is deliberate, since such maps preserve the intrinsic rank of the moment pair (see Proposition~\ref{prop:rankreduction-carrying}).

\subsection{Radial cut profile}
For brevity, we use the following terminology.
\begin{definition}[Radial cut profile]\label{def:profile}
A \emph{radial cut profile}, or \emph{radial profile} for short, is a finite-valued nonincreasing function
\[
   \rho:[0,1]\longrightarrow[0,\infty)
\]
that is left-continuous on \((0,1]\) and satisfies
\[
   \rho(0)=\lim_{\alpha\downarrow0}\rho(\alpha).
\]
\end{definition}
\subsection{Characterization by affine consistency}
The following is the main characterization result of this paper. It shows that affine consistency fixes the geometry of every cut, leaving only its radial profile undetermined.

\begin{theorem}[Characterization by affine consistency]\label{thm:characterization}
A family of maps \(\Phi_d:\mathcal M_d\to\F_d\), \(d\ge1\), is affine-consistent if and only if there exists a unique family of radial cut profiles \((\rho_p)_{p\ge1}\) such that  \(\Phi_d(\bmu,0)\) is the crisp singleton \(\{\bmu\}\) and, whenever \(p=\rank\bSigma\ge1\), for all  $0\le\alpha\le1$,
\begin{equation}\label{eq:characterization}
   C_\alpha\bigl(\Phi_d(\bmu,\bSigma)\bigr)
      =
   \bmu+
   \rho_p(\alpha)\bSigma^{1/2}
   \mathbb B_{\Img\bSigma},
\end{equation}
equivalently,
\[
   C_\alpha\bigl(\Phi_d(\bmu,\bSigma)\bigr)
      =
   \bigl\{
      \bx\in\Sspace(\bmu,\bSigma):
      t_{\bmu,\bSigma}(\bx)\le\rho_p(\alpha)
   \bigr\}.
\]
\end{theorem}
From the above,  the profile is indexed by the intrinsic rank \(p\), that is, the covariance rank, alone. It does not depend on the ambient dimension, the orientation of the state space, or the eigenvalues of \(\bSigma\).

\begin{proof}
For ease of reading, we present the necessity part of the proof in four steps.

\smallskip
\noindent
\emph{Step 1: Zero-variance directions.}
Fix \(d\), \(1\le p\le d\), and \(\alpha\in[0,1]\), and consider the canonical rank-\(p\) covariance matrix
\[
   \bJ_{d,p}:=\operatorname{diag}(\bI_p,0_{d-p}).
\]
Put
\[
   K_\alpha:=C_\alpha\bigl(\Phi_d(0,\bJ_{d,p})\bigr).
\]
For every \(s>0\), the matrix
\[
   \bS_s:=\operatorname{diag}(\bI_p,s\bI_{d-p})
\]
satisfies \(\bS_s\bJ_{d,p}\bS_s^\top=\bJ_{d,p}\). Coordinate neutrality and \eqref{eq:cutimage} therefore give
\[
   \bS_sK_\alpha=K_\alpha
   \qquad\text{for every }s>0.
\]
If \(K_\alpha\) contained a point \((\bu,\bv)\) with \(\bv\ne0\), it would also contain \((\bu,s\bv)\) for every \(s>0\), contradicting compactness. Hence
\begin{equation}\label{eq:canonicalstate}
   K_\alpha\subseteq\R^p\times\{0\}.
\end{equation}

\medskip
\noindent
\emph{Step 2: Rotational symmetry.}
Let \(\Orth(p):=\{\bQ\in\R^{p\times p}:\bQ^\top\bQ=\bI_p\}\) denote the group of orthogonal \(p\times p\) matrices, that is, of rotations and reflections of \(\R^p\). For every \(\bQ\in\Orth(p)\), the transformation
\[
   \widetilde\bQ:=\operatorname{diag}(\bQ,\bI_{d-p})
\]
also leaves the moment pair \((0,\bJ_{d,p})\) unchanged. Hence \( \widetilde\bQ K_\alpha=K_\alpha\). Viewed in its effective \(p\)-dimensional space, \(K_\alpha\) is therefore a nonempty compact convex set invariant under every rotation and reflection about the origin.

Let
\[
   r_{d,p}(\alpha):=\max\{\|\bx\|:\bx\in K_\alpha\}.
\]
If \(\bx\in K_\alpha\) has norm \(r_{d,p}(\alpha)\), rotational invariance puts the whole sphere of that radius in \(K_\alpha\). Convexity then puts its convex hull, the corresponding ball, in \(K_\alpha\), and maximality of the radius gives the reverse inclusion. Thus
\begin{equation}\label{eq:canonicalball}
   K_\alpha
      =
   r_{d,p}(\alpha)\mathbb B_p\times\{0\}.
\end{equation}
The geometric mechanism is therefore
\[
   \text{rotational invariance}
   +\text{convexity}
   +\text{compactness}
   \quad\Longrightarrow\quad
   \text{Euclidean ball}.
\]

\medskip
\noindent
\emph{Step 3: General covariance by affine transport.}
Let now \(\bSigma\succeq0\) have rank \(p\). There exists an invertible matrix \(\bC\in\GL(d)\) such that
\[
   \bSigma=\bC\bJ_{d,p}\bC^\top.
\]
Applying coordinate neutrality to the affine map \(\bx\mapsto \bC\bx+\bmu\) gives
\[
   C_\alpha\bigl(\Phi_d(\bmu,\bSigma)\bigr)
      =
   \bmu+r_{d,p}(\alpha)\bC\bigl(\mathbb B_p\times\{0\}\bigr).
\]
The support function of the centered set \(\bC(\mathbb B_p\times\{0\})\) is
\[
   \bu
   \longmapsto
   \sqrt{\bu^\top \bC\bJ_{d,p}\bC^\top \bu}
   =
   \sqrt{\bu^\top\bSigma \bu}.
\]
By \eqref{eq:supportaffine}, and since \(\bSigma^{1/2}\bu\in\Img\bSigma\), this is also the support function of \(\bSigma^{1/2}\mathbb B_{\Img\bSigma}\). Since a compact convex set is uniquely characterized by its support function,
\[
   \bC\bigl(\mathbb B_p\times\{0\}\bigr)
      =
   \bSigma^{1/2}\mathbb B_{\Img\bSigma}.
\]
Hence
\begin{equation}\label{eq:rankprofile}
   C_\alpha\bigl(\Phi_d(\bmu,\bSigma)\bigr)
      =
   \bmu+r_{d,p}(\alpha)\bSigma^{1/2}\mathbb B_{\Img\bSigma}.
\end{equation}

\medskip
\noindent
\emph{Step 4: Ambient neutrality and intrinsic rank.}
Nestedness of the cuts implies that \(\alpha\mapsto r_{d,p}(\alpha)\) is nonincreasing, while \eqref{eq:cutregularity} gives left-continuity on \((0,1]\) and
\[
   r_{d,p}(0)=\lim_{\alpha\downarrow0}r_{d,p}(\alpha).
\]
Thus \(r_{d,p}\) is a radial cut profile.

Ambient neutrality removes the remaining dependence on the ambient dimension. By \eqref{eq:embedding} with \((p,d)\) in place of \((d,m)\),
\[
   \Phi_d(0,\bJ_{d,p})=\iota_{p,d}\cdot\Phi_p(0,\bI_p).
\]
Taking cuts with \eqref{eq:cutimage} and using \eqref{eq:canonicalball} in both dimensions gives \(r_{d,p}(\alpha)\mathbb B_p\times\{0\}=r_{p,p}(\alpha)\mathbb B_p\times\{0\}\), so
\[
   r_{d,p}=r_{p,p}.
\]
Writing \(\rho_p:=r_{p,p}\) gives \eqref{eq:characterization}.

If \(p=0\), the moment pair \((0,0)\) is unchanged by every scaling \(\bx\mapsto s\bx\), \(s>0\). As in Step 1, every cut \(K\) of \(\Phi_d(0,0)\) then satisfies \(sK=K\) for all \(s>0\), so boundedness and nonemptiness force \(K=\{0\}\), and translation gives the crisp singleton \(\{\bmu\}\).

\medskip
\noindent
\emph{Sufficiency.}
Conversely, let \((\rho_p)_{p\ge1}\) be any sequence of radial cut profiles and define the cuts by \eqref{eq:characterization}, with the singleton convention when \(\bSigma=0\). Their nesting and regularity give, by the converse following \eqref{eq:cutregularity}, a fuzzy vector in \(\F_d\).

Let
\[
   T(\bx)=\bA\bx+\bb:\R^d\longrightarrow\R^m
\]
be injective. If \(\bSigma=0\), then \(\bA\bSigma \bA^\top=0\) and \(T\cdot\{\bmu\}=\{\bA\bmu+\bb\}\), so \eqref{eq:injective} holds trivially. Otherwise \(\rank(\bA\bSigma \bA^\top)=\rank\bSigma=p\ge1\). The support function of the image of the centered cut in \eqref{eq:characterization} is
\[
   \bu
   \longmapsto
   \rho_p(\alpha)\sqrt{\bu^\top \bA\bSigma \bA^\top \bu}.
\]
This is also the support function of
\[
   \rho_p(\alpha)(\bA\bSigma \bA^\top)^{1/2}
   \mathbb B_{\Img(\bA\bSigma \bA^\top)}.
\]
The corresponding cuts therefore coincide, and \eqref{eq:injective} follows. Hence the family is affine-consistent by Lemma~\ref{lem:injective}.

In either direction, the profiles are unique, since the cuts of \(\Phi_p(0,\bI_p)\) are \(\rho_p(\alpha)\mathbb B_p\).
\end{proof}

\paragraph{Notation.}
For a radial cut profile \(\rho\), we write
\[
   \Ell(\bmu,\bSigma,\rho)
\]
for the fuzzy vector with cuts
\begin{equation}\label{eq:ellipsoidalcuts}
   C_\alpha
      =
   \bmu+\rho(\alpha)\bSigma^{1/2}\mathbb B_{\Img\bSigma},
\end{equation}
with the crisp singleton convention when \(\bSigma=0\). In particular, the crisp Mahalanobis ellipsoid \(E_r(\bmu,\bSigma)\) is \(\Ell(\bmu,\bSigma,\rho)\) with \(\rho\equiv r\). Equivalently,
\begin{equation}\label{eq:ellipsoidalmembership}
   \pi_{\Ell}(\bx)
      =
   \sup\{\alpha\in[0,1]:t_{\bmu,\bSigma}(\bx)\le\rho(\alpha)\},
   \qquad \sup\emptyset=0.
\end{equation}
Theorem~\ref{thm:characterization} says precisely that
\[
   \Phi_d(\bmu,\bSigma)=\Ell(\bmu,\bSigma,\rho_p)
\]
whenever \(p=\rank\bSigma\ge1\). At this stage no carrying requirement has been imposed. In particular, \(\rho_p\equiv0\) is compatible with affine consistency. The theorem determines the form of the cuts, not yet how large they must be.

\subsection{Affine images and linear combinations}\label{sec:interaction}
For any affine map \(T(\bx)=\bA\bx+\bb\) from \(\R^d\) to \(\R^{d'}\) and any fixed radial cut profile \(\rho\), the support-function identity \eqref{eq:supportaffine} gives
\begin{equation}\label{eq:ellimage}
 T\cdot\Ell(\bmu,\bSigma,\rho)
   =
 \Ell(\bA\bmu+\bb,\bA\bSigma \bA^\top,\rho).
\end{equation}
This holds for every \(\bA\), including rank-reducing maps, and the profile on the right is \(\rho\) itself. For an affine-consistent family, let \(p=\rank\bSigma\) and \(p'=\rank(\bA\bSigma\bA^\top)\ge1\). Then \eqref{eq:ellimage} gives \(T\cdot\Phi_d(\bmu,\bSigma)=\Ell(\bA\bmu+\bb,\bA\bSigma\bA^\top,\rho_p)\), whereas \(\Phi_{d'}(\bA\bmu+\bb,\bA\bSigma\bA^\top)=\Ell(\bA\bmu+\bb,\bA\bSigma\bA^\top,\rho_{p'})\). The two coincide if and only if \(\rho_p=\rho_{p'}\). Affine consistency concerns injective maps only and imposes no such relation, which is why the profiles remain indexed by the intrinsic rank.

In particular, for a joint ellipsoidal fuzzy vector \(\xi=\Ell(\bmu,\bSigma,\rho)\) and any \(\ba\in\R^d\),
\begin{equation}\label{eq:linearcombo}
 C_\alpha(\ba^\top\cdot\xi)
   =
 [\ba^\top\bmu\pm\rho(\alpha)\sqrt{\ba^\top\bSigma \ba}],
\end{equation}
where the bracket with \(\pm\) denotes the two endpoints. Thus the ellipsoidal joint aggregates marginal spreads and interaction through the covariance quadratic form. The minimum-based non-interactive joint of the same one-dimensional marginals has rectangular cuts and gives the half-width \(\rho(\alpha)\sum_i|a_i|\sqrt{\Sigma_{ii}}\), which is at least as large by the triangle inequality. The difference persists for uncorrelated coordinates. For \(\bmu=0\), \(\bSigma=\bI_2\) and \(\ba=(1,1)^\top\), the two half-widths are \(\rho(\alpha)\sqrt2\) and \(2\rho(\alpha)\). The distinction therefore reflects the joint possibility structure, not only nonzero cross-covariances \cite{DuboisPrade1981,FullerMajlender2004,CarlssonFullerMajlender2005}.

\subsection{Role of the assumptions}\label{sec:role-assumptions}

Each assumption of Theorem~\ref{thm:characterization} is needed. The examples below drop one of them at a time, coordinate neutrality~\ref{prop:coordinateneut} or ambient neutrality~\ref{prop:ambientneut} from Definition~\ref{def:affine}, or convexity or compactness of the cuts, while keeping the others. In each case a conclusion of the theorem fails.

\smallskip
\paragraph{Ambient neutrality without coordinate neutrality}
Consider the crisp coordinate box
\[
   K(\bmu,\bSigma)
      :=
   \prod_{i=1}^d
   \bigl[
      \mu_i-\sqrt{\Sigma_{ii}},
      \mu_i+\sqrt{\Sigma_{ii}}
   \bigr].
\]
This construction depends only on the moment pair and is unaffected by adding deterministic zero coordinates, but it is not coordinate-neutral. For \((\bmu,\bSigma)=(0,\bI_2)\) it gives the square \([-1,1]^2\). A rotation \(\bQ\) by \(\pi/4\) leaves the moment pair unchanged, so the construction returns the same square, whereas coordinate neutrality would require the rotated square \(\bQ[-1,1]^2\).

\smallskip
\paragraph{Coordinate neutrality without ambient neutrality}
Within each fixed dimension \(d\), the cuts
\[
   C_\alpha^{(d)}(\bmu,\bSigma)
      =
   \bmu+\sqrt d\,(1-\alpha)\bSigma^{1/2}\mathbb B_{\Img\bSigma}
\]
define a coordinate-neutral construction. After embedding the same moment pair into a higher-dimensional space, however, the radial factor changes from \(\sqrt d\) to \(\sqrt m\). Ambient neutrality rules out exactly such dependence on coordinates that carry no additional uncertainty.

\smallskip
\paragraph{Affine consistency without convexity}
Assigning to each moment pair the crisp Mahalanobis sphere
\[
   \{\bx\in\Sspace(\bmu,\bSigma):t_{\bmu,\bSigma}(\bx)=1\}
\]
gives compact cuts and satisfies both conditions of Definition~\ref{def:affine} on the larger class of compact fuzzy sets, with the singleton convention at rank zero. More generally, rotationally invariant shells become possible. Thus symmetry alone does not imply a ball. Convexity fills the sphere into its convex hull.

\smallskip
\paragraph{Affine consistency without compactness}
Assigning to each moment pair its entire state space \(\Sspace(\bmu,\bSigma)\), viewed as a crisp fuzzy set, satisfies both conditions of Definition~\ref{def:affine} on the larger class of closed convex fuzzy sets, with convex cuts but no finite radial profile when \(p>0\). Compactness supplies a finite maximal radius.

The characterization is therefore conditional on a visible set of assumptions. Coordinate and ambient neutrality determine the invariant geometry, while convexity and compactness determine the admissible form of the cuts. Nestedness makes the radius nonincreasing in \(\alpha\), while \eqref{eq:cutregularity} gives its one-sided continuity and its value at \(\alpha=0\). No differentiability of the radial profile or distributional assumption is involved.

\smallskip
The remaining freedom is the radial profile. The next section asks what changes when the fuzzy vector is additionally required to carry the prescribed moments.

\section{Characterizing carrying through the radial profile}\label{sec:carrying}
We now connect Theorem~\ref{thm:characterization}, which is purely geometric, to the probabilistic meaning of the moments, and ask when the carrying requirement of Definition~\ref{def:carrying} is satisfied.

Throughout this section, $p=\rank\bSigma\ge1$. We call a selection $\bX$ in the sense of Definition~\ref{def:selection} \emph{realizing} if $\E\bX=\bmu$ and $\Cov\bX=\bSigma$. Since such an $\bX$ lies in the state space almost surely, the moments imply the trace identity
\begin{equation}\label{eq:trace}
 \E\,t(\bX)^2=\tr\!\left(\bSigma^+\E[(\bX-\bmu)(\bX-\bmu)^\top]\right)
 =\tr(\bSigma^+\bSigma)=p.
\end{equation}
In an ellipsoidal fuzzy vector $\Ell(\bmu,\bSigma,\rho)$, the cuts \eqref{eq:ellipsoidalcuts} are the sets $\{\bx\in\Sspace(\bmu,\bSigma):t(\bx)\le\rho(\alpha)\}$, so every selection satisfies
\begin{equation}\label{eq:radialbound}
t(\bX)\le\rho(U)
\end{equation}
almost surely. This yields an exact criterion as our second main result.

\begin{theorem}[Carrying and normalization]\label{thm:carrying}
Let $\xi=\Ell(\bmu,\bSigma,\rho)$, where $\rho$ is a radial cut profile and $p\ge1$.
\begin{enumerate}[label=(\roman*)]
\item $\xi$ carries $(\bmu,\bSigma)$ if and only if
\begin{equation}\label{eq:carrybound}
 \int_0^1\rho(\alpha)^2\dd\alpha\ge p.
\end{equation}
\item If equality holds in \eqref{eq:carrybound}, which we call \emph{normalization}, the canonical selection
\begin{equation}\label{eq:canonical}
 \bX=\bmu+\rho(U)\bSigma^{1/2}\bTheta,\quad
 U\sim\Unif(0,1),\quad \bTheta\sim\Unif(\Sphere_{\Img\bSigma}),
\end{equation}
with $U$ and $\bTheta$ independent, realizes the moments. Every realizing selection then satisfies $t(\bX)=\rho(U)$ almost surely.
\item Within the class $\{\Ell(\bmu,\bSigma,\rho):\rho\text{ a radial cut profile}\}$ for fixed $(\bmu,\bSigma)$, ordered by inclusion of all cuts, equality in \eqref{eq:carrybound} characterizes the minimal carriers. None of their cuts can be reduced, while keeping the others no larger and remaining in this class, without losing carrying.
\end{enumerate}
\end{theorem}
\begin{proof}
Necessity follows from \eqref{eq:trace} and \eqref{eq:radialbound}, which give $\int_0^1\rho^2\ge\E\,t(\bX)^2=p$. For sufficiency write $c:=\int_0^1\rho^2\ge p$ and $\kappa:=\sqrt{p/c}\le1$. Let $U$ and $\bTheta$ be independent as in \eqref{eq:canonical}, and set $\bX=\bmu+\kappa\rho(U)\bSigma^{1/2}\bTheta$. It is a selection because $t(\bX)=\kappa\rho(U)\le\rho(U)$. Rotational symmetry yields $\E\bTheta=0$ and $\E(\bTheta\bTheta^\top)=p^{-1}\bPi_{\bSigma}$. It follows that $\E \bX=\bmu$ and $\Cov \bX=\kappa^2c\bSigma/p=\bSigma$.

When $c=p$, we have $\kappa=1$ and \eqref{eq:canonical} applies. For any realizing selection, $\rho(U)^2-t(\bX)^2\ge0$ and its expectation is zero by \eqref{eq:trace}. Thus $t(\bX)=\rho(U)$ almost surely. Finally, if $\rho'\le\rho$ is strictly smaller at one level, left-continuity (or the limit condition at zero) gives strict inequality on an interval of levels, so $\int_0^1(\rho')^2<\int_0^1\rho^2$. Conversely, if $c>p$, multiplying $\rho$ by $\sqrt{p/c}<1$ strictly reduces its nonzero cuts and preserves carrying. This proves the minimality claim.
\end{proof}
Equality forces saturation of the available radius. On every nondegenerate random cut, a realizing selection is on its boundary relative to the state space. When $\rho(U)=0$, the cut is simply $\{\bmu\}$. The canonical selection \eqref{eq:canonical} is an elliptically contoured law with radial part $\rho(U)$ \cite{CambanisHuangSimons1981}, supported on $\{\bx\in\Sspace(\bmu,\bSigma):t(\bx)\in\overline{\rho((0,1))}\}$. Its covariance is $\int_0^1\rho(\alpha)^2\dd\alpha\,/p \cdot \bSigma$, so normalization is exactly the condition under which it equals $\bSigma$.

For a family that is affine-consistent according to Definition~\ref{def:affine}, the condition reads
\begin{equation}\label{eq:familycarry}
 \int_0^1\rho_p(\alpha)^2\dd\alpha\ge p\quad\text{for every }p\ge1.
\end{equation}
It also implies $\rho_p(0)\ge\sqrt p$. A family is \emph{normalized} at rank $p$ if $\rho_p$ is normalized in the sense of Theorem~\ref{thm:carrying}(ii).

Normalization selects minimal carriers within the fixed-shape Mahalanobis class, not one uniquely smallest carrier. For example, for any $\beta>0$ the profiles
\begin{equation}\label{eq:betaprofile}
 \rho_{\beta,p}(\alpha)=\sqrt{\frac{p(\beta+2)}{\beta}}(1-\alpha)^{1/\beta}
\end{equation}
are normalized. The triangular membership $\pi(\bx)=(1-t(\bx)/\sqrt{3p})_+$ corresponds to $\beta=1$, and the parabolic membership $\pi(\bx)=(1-t(\bx)^2/(2p))_+$ to $\beta=2$.
We illustrate Theorem~\ref{thm:carrying} for a concrete numerical rank-two covariance, constructing two normalized profiles with explicit realizing selections. The example also shows that normalization does not determine the radial profile uniquely.

\begin{example}[Fixed geometry, non-unique minimal profiles]\label{ex:bivariate}
We work in the plane, $d=2$, with the moment pair
\begin{equation}\label{eq:examplemoments}
 \bmu=0,\qquad \bSigma=\begin{pmatrix}1&3/5\\3/5&1\end{pmatrix}.
\end{equation}
Both variances equal $1$ and the correlation is $3/5$. Since $\det\bSigma=1-9/25=16/25>0$, the covariance is nonsingular. Hence $p=2$, and the state space is all of $\R^2$.

\smallskip
\noindent
\emph{Step 1: Mahalanobis radius.}
Because $\bSigma$ is invertible, $\bSigma^+=\bSigma^{-1}$, and
\[
 \bSigma^{-1}=\frac{25}{16}\begin{pmatrix}1&-3/5\\-3/5&1\end{pmatrix}.
\]
Therefore
\[
 t(\bx)^2=\bx^\top\bSigma^{-1}\bx=\frac{25}{16}\left(x_1^2-\frac65x_1x_2+x_2^2\right).
\]

\smallskip
\noindent
\emph{Step 2: Shape of the cuts.}
The eigenvectors of $\bSigma$ are $(1,1)^\top/\sqrt2$ and $(1,-1)^\top/\sqrt2$, with eigenvalues $1+3/5=8/5$ and $1-3/5=2/5$. Its square root is
\[
 \bB:=\bSigma^{1/2}=\frac1{\sqrt{10}}\begin{pmatrix}3&1\\1&3\end{pmatrix},
\]
which has the same eigenvectors and the eigenvalues $\sqrt{8/5}$ and $\sqrt{2/5}$. A direct multiplication confirms $\bB^2=\bSigma$. By Theorem~\ref{thm:characterization}, every cut is
\[
 C_\alpha=\rho(\alpha)\bB\,\mathbb B_2,
\]
the unit disc stretched by $\bB$ and then scaled by $\rho(\alpha)$. It is an ellipse with principal axes along $(1,1)$ and $(1,-1)$ and semiaxis lengths $\rho(\alpha)\sqrt{8/5}$ and $\rho(\alpha)\sqrt{2/5}$. The orientation and the axis ratio $2:1$ are therefore fixed by $\bSigma$.

Only the common factor $\rho(\alpha)$ remains free.

\smallskip
\noindent
\emph{Step 3: Two normalized profiles.}
Consider the two profiles
\begin{equation}\label{eq:exampleprofiles}
 \rho_{\mathrm{tri}}(\alpha)=\sqrt6(1-\alpha),\qquad
 \rho_{\mathrm{par}}(\alpha)=2\sqrt{1-\alpha}.
\end{equation}
Both are continuous and nonincreasing and vanish at $\alpha=1$, so they are radial cut profiles in the sense of Definition~\ref{def:profile}. Normalization requires $\int_0^1\rho^2\dd\alpha=p=2$, and indeed
\begin{align*}
 \int_0^1\rho_{\mathrm{tri}}(\alpha)^2\dd\alpha&=6\int_0^1(1-\alpha)^2\dd\alpha=6\cdot\frac13=2,\\
 \int_0^1\rho_{\mathrm{par}}(\alpha)^2\dd\alpha&=4\int_0^1(1-\alpha)\dd\alpha=4\cdot\frac12=2.
\end{align*}
By Theorem~\ref{thm:carrying}, both profiles give minimal carriers of \eqref{eq:examplemoments}.

\smallskip
\noindent
\emph{Step 4: Membership functions.}
For $\alpha\in(0,1]$, a point $\bx$ has membership $\alpha$ exactly when it lies on the boundary of the cut $C_\alpha$, that is, when $t(\bx)=\rho(\alpha)$. Solving this equation for $\alpha$ gives
\[
 \pi_{\mathrm{tri}}(\bx)=\bigl(1-t(\bx)/\sqrt6\bigr)_+,
 \qquad \pi_{\mathrm{par}}(\bx)=\bigl(1-t(\bx)^2/4\bigr)_+.
\]
The first membership decreases linearly in the Mahalanobis radius and the second quadratically. They are the cases $\beta=1$ and $\beta=2$ of \eqref{eq:betaprofile} with $p=2$.

\smallskip
\noindent
\emph{Step 5: Comparing the two carriers.}
Setting $\sqrt6(1-\alpha)=2\sqrt{1-\alpha}$ and dividing by $\sqrt{1-\alpha}$ gives $\sqrt{1-\alpha}=2/\sqrt6$, that is, $\alpha=1/3$. At this level both cuts have radius $2\sqrt6/3$. At $\alpha=0$ the triangular profile is larger, since $\sqrt6>2$. Hence for $0\le\alpha<1/3$ the triangular profile gives larger cuts, and for $1/3<\alpha<1$ it gives smaller cuts. Both vanish at $\alpha=1$. Thus neither fuzzy vector contains the other cut by cut. Minimality is not uniqueness, and there is no least carrier in this fixed-shape class, because such a carrier would have to lie below both incomparable minimal elements.

\smallskip
\noindent
\emph{Step 6: An explicit realizing selection.}
The construction does not require choosing a probability law in advance. Let $U\sim\Unif(0,1)$ and $V\sim\Unif(0,2\pi)$ be independent, and set
\[
 \bX=\rho(U)\bB\bTheta,
 \qquad
 \bTheta=\begin{pmatrix}\cos V\\\sin V\end{pmatrix},
\]
for either profile $\rho$. The vector $\bTheta$ is uniform on the unit circle. Since $\E\cos V=\E\sin V=\E\cos V\sin V=0$ and $\E\cos^2V=\E\sin^2V=1/2$, we have $\E\bTheta=0$ and $\E\bTheta\bTheta^\top=\frac12\bI_2$. Using independence and normalization,
\[
 \E\bX=0,
 \qquad
 \Cov\bX=\E\bigl[\rho(U)^2\bigr]\,\bB\Bigl(\tfrac12\bI_2\Bigr)\bB=2\cdot\tfrac12\,\bB^2=\bSigma.
\]
Moreover $t(\bX)=\|\bB^{-1}\bX\|=\rho(U)\|\bTheta\|=\rho(U)$, so $\bX$ lies on the boundary of its random cut $C_U$, as Theorem~\ref{thm:carrying}(ii) requires. Hence $\bX$ is a realizing selection for either profile. It is the canonical selection \eqref{eq:canonical} with $p=2$.
\end{example}

\begin{remark}[\cite{FullerMajlender2003} variance]%
For $p=d=1$ and $\bSigma=\sigma^2>0$, normalization has a possibilistic reading. The fuzzy number $\Ell(\mu,\sigma^2,\rho)$ has cuts $[\mu-\rho(\alpha)\sigma,\mu+\rho(\alpha)\sigma]$, so its $f$-weighted possibilistic variance in the sense of Full\'er and Majlender \cite{FullerMajlender2003} is $\sigma^2\int_0^1\rho(\alpha)^2f(\alpha)\dd\alpha$. With the constant weighting $f\equiv1$, normalization says exactly that this variance equals $\sigma^2$. Other weightings give other scales. For instance, $f(\alpha)=2\alpha$ gives $\sigma^2/2$ for the normalized triangular profile $\rho_{1,1}$.
\end{remark}

\begin{proposition}[Rank reduction and universal carrying]
\label{prop:rankreduction-carrying}
Let $\Phi=(\Phi_d)_{d\ge1}$ be an affine-consistent family
as in Theorem~\ref{thm:characterization}. Suppose, in
addition, that $\Phi$ commutes with all affine maps,
including rank-reducing ones. Then:
\begin{enumerate}[label=(\roman*)]
\item All positive-rank radial profiles coincide:
      $\rho_p=\rho$ for every $p\ge1$.
\item No such family can carry every
      moment pair of every rank.
\item Such a family cannot be normalized at two distinct
      positive covariance ranks.
\end{enumerate}
\end{proposition}

\begin{proof}
For $p\ge2$, consider the canonical moment pair
$(0,\bI_p)$ and the projection
$\bA:\R^p\to\R$ onto the first coordinate.
By Theorem~\ref{thm:characterization}, the assigned
fuzzy vector has radial profile $\rho_p$.
Its image under $\bA$ retains this profile by
\eqref{eq:ellimage}, whereas the assignment to the
projected moment pair $(0,1)$ has profile $\rho_1$.
As noted after \eqref{eq:ellimage}, full affine compatibility therefore implies
$\rho_p=\rho_1$ for every $p\ge2$, proving (i).

Write $\rho$ for the common profile and put
\[
   c:=\int_0^1\rho(\alpha)^2\dd\alpha.
\]
Since $\rho$ is bounded, $c<\infty$.
By Theorem~\ref{thm:carrying}(i), carrying a moment
pair of rank $p$ requires $c\ge p$.
This cannot hold for every positive integer $p$,
which proves (ii).

Finally, normalization at rank $p$ means $c=p$ by definition.
The same profile cannot satisfy this equality at
two distinct positive ranks, proving (iii).
\end{proof}

The proposition shows that restricting affine consistency
to injective maps is substantive rather than merely
technical. Full affine compatibility eliminates the
freedom to choose rank-dependent profiles, but thereby
conflicts with universal carrying across all ranks.
For normalized carriers, the incompatibility already
arises between ranks one and two.

The marginals in the comparison of Section~\ref{sec:interaction} are the projections \(\be_i^\top\cdot\xi=\Ell(\mu_i,\Sigma_{ii},\rho)\) of the joint, which keep the profile \(\rho\). For coordinates with positive variance, they differ from the fuzzy numbers \(\Phi_1(\mu_i,\Sigma_{ii})=\Ell(\mu_i,\Sigma_{ii},\rho_1)\) that an affine-consistent family assigns to each coordinate separately whenever \(\rho\ne\rho_1\). Zero-variance coordinates remain crisp singletons. For a family normalized at ranks \(1\) and \(p\ge2\), necessarily \(\rho_p\ne\rho_1\). For the profiles \eqref{eq:betaprofile}, for instance, \(\rho_{\beta,p}=\sqrt p\,\rho_{\beta,1}\), so the half-widths of the marginals of a rank-\(p\) joint are wider than the separately constructed ones by the factor \(\sqrt p\), and the comparison concerns the former.

\section{Sparsest realizations and symmetry completion}\label{sec:simplex}
The canonical selections of the previous section are rotationally symmetric. This section now turns to the opposite extreme, discrete realizations of the moments by as few points (atoms) as possible, and shows how their geometry relates to the ellipsoidal carriers.

\smallskip

An atomic law with mean $\bmu$ and covariance $\bSigma$ is a discrete realization of the moment pair. Since its support must affinely span the $p$-dimensional state space, at least $p+1$ atoms are needed. We call a realization with exactly $p+1$ atoms \emph{sparsest}. The covariance-scenario construction \cite{FastEmpirical} obtains exactly $p+1$ equally weighted points $\bzeta_1,\ldots,\bzeta_{p+1}$ with a single Householder reflection, and its details are not repeated here. What matters for our purpose is the geometry of every equally weighted sparsest realization, independently of how its points are computed. Equal weights are a requirement of this particular realization, not a consequence of having only $p+1$ atoms.\footnote{
The moment pair $(\bmu,\bSigma)$ need not arise from an atomic law. It may be the mean and covariance of any probability law, continuous or discrete, or it may be specified directly without reference to a law. Whatever its origin, it can always be represented by an atomic law. A discrete realization is thus a finitely atomic representing measure in the terminology of the truncated moment problem \cite{CurtoFialkow2005}. It represents the moments but is not assumed to be the law that generated them.}

\smallskip
Let $\mathcal Z=\{\bzeta_1,\ldots,\bzeta_{p+1}\}$ and $\mathcal T:=\conv\mathcal Z$. In Mahalanobis coordinates write $\bz_j:=\bSigma^{+/2}(\bzeta_j-\bmu)\in\Img\bSigma$. Let $\Orth(\Img\bSigma)$ denote the orthogonal group of $\Img\bSigma$, that is, the set of linear maps $\bQ:\Img\bSigma\to\Img\bSigma$ with $\norm{\bQ\bz}=\norm{\bz}$ for all $\bz\in\Img\bSigma$. In an orthonormal basis of $\Img\bSigma$ these maps are represented by the matrices in $\Orth(p)$. Every $\bQ\in\Orth(\Img\bSigma)$ induces a map on the state space,
\[
 R_{\bQ}(\bmu+\bSigma^{1/2}\bz)=\bmu+\bSigma^{1/2}\bQ\bz.
\]
Such Mahalanobis rotations preserve the moment pair and extend to invertible affine transformations of the ambient space by acting as the identity on $\ker\bSigma$. The term includes reflections. We offer the following result.

\begin{theorem}[Regular simplex and rotational completion]\label{thm:simplex}
Suppose $p\ge1$ and $\bzeta_1,\ldots,\bzeta_{p+1}$ are equally weighted points with moments $(\bmu,\bSigma)$. Then the following hold.
\begin{enumerate}[label=(\roman*)]
\item They are affinely independent, and their Mahalanobis coordinates form a regular $p$-simplex with
\begin{equation}\label{eq:gram}
 \norm{\bz_j}^2=p,\qquad \langle \bz_i,\bz_j\rangle=-1\quad(i\ne j).
\end{equation}
Any two labeled such configurations differ by a unique Mahalanobis rotation.
\item The union of the rotational images of the simplex is
\begin{equation}\label{eq:completion}
 \bigcup_{\bQ\in\Orth(\Img\bSigma)}R_{\bQ}(\mathcal T)
 =\bmu+\sqrt p\,\bSigma^{1/2}\mathbb B_{\Img\bSigma}.
\end{equation}
The union is already convex and is the smallest crisp Mahalanobis ellipsoid that carries $(\bmu,\bSigma)$. It is also the unique ellipsoid in the state space of minimal volume that carries $(\bmu,\bSigma)$.
\end{enumerate}
\end{theorem}
\begin{proof}
The $p+1$ points must affinely span the $p$-dimensional state space, so they are affinely independent. Choose an orthonormal basis of $\Img\bSigma$ and form the $(p+1)\times p$ matrix $\bZ$ with rows $\bz_j^\top$. The moment equations become $\bZ^\top\mathbf1=0$ and $\bZ^\top \bZ=(p+1)\bI_p$. The Gram matrix $\bG=\bZ\bZ^\top$ is $(p+1)$ times an orthogonal projector with kernel $\operatorname{span}\{\mathbf1\}$, hence $\bG=(p+1)\bI_{p+1}-\mathbf1\mathbf1^\top$. Its diagonal and off-diagonal entries give \eqref{eq:gram}. Two labeled configurations with the same Gram matrix are related by a unique orthogonal map on their common spanning space.

In Mahalanobis coordinates $\mathcal T$ lies in the ball of radius $\sqrt p$ and contains both the origin and a vertex at that radius. It therefore contains points of every radius between $0$ and $\sqrt p$. The orthogonal group acts transitively on each sphere of fixed radius, so its orbit covers the entire ball and nothing outside it, proving \eqref{eq:completion}. The assertion on the smallest crisp Mahalanobis ellipsoid follows from Theorem~\ref{thm:carrying} with constant profile $\rho\equiv c$, since the crisp ellipsoid $\{t\le c\}$ carries the moments precisely when $c\ge\sqrt p$. For the volume statement, work in Mahalanobis coordinates, where the moments become $(0,\bI_p)$ and all volumes in the state space are multiplied by the same factor. A carrier must have full dimension $p$ by \eqref{eq:affinesupport}, so consider $\{\bz:(\bz-\bc)^\top\bA^{-1}(\bz-\bc)\le1\}$ with $\bA\succ0$. A realizing law $\bZ$ satisfies $1\ge\E(\bZ-\bc)^\top\bA^{-1}(\bZ-\bc)=\tr\bA^{-1}+\bc^\top\bA^{-1}\bc$. By the inequality of arithmetic and geometric means, $\tr\bA^{-1}\ge p(\det\bA)^{-1/p}$, so $\det\bA\ge p^p$, with equality only for $\bc=0$ and $\bA=p\bI_p$. The volume is proportional to $(\det\bA)^{1/2}$, which proves the claim.
\end{proof}
\begin{remark}There is a probabilistic counterpart. Let $\bQ$ be Haar-uniform on $\Orth(\Img\bSigma)$, that is, drawn from its unique invariant probability law \cite{Haar1933,Meckes2019}. If $\bQ$ is independent of a uniform vertex index $N$, then $R_{\bQ}(\bzeta_{N})$ has the same law as $\bmu+\sqrt p\bSigma^{1/2}\bTheta$, with $\bTheta\sim\Unif(\Sphere_{\Img\bSigma})$. This is averaging of probability laws through Haar randomization. The geometric operation in \eqref{eq:completion}, by contrast, is the union of rotated sets. The randomized vertex law is the canonical selection of the crisp normalized ellipsoid. The ellipsoid in \eqref{eq:completion} is also the minimum-volume ellipsoid in the state space that contains $\mathcal T$, by the duality between minimum-volume enclosing ellipsoids and D-optimal designs \cite{Titterington1975,Todd2016}, since all vertices of $\mathcal T$ have the same Mahalanobis radius $\sqrt p$. Theorem~\ref{thm:simplex}(ii) is a statement of a different kind, because it concerns ellipsoids that admit a probability law with the prescribed moments rather than ellipsoids that contain given points.
\end{remark}
The same construction extends from the crisp simplex to fuzzy vectors. Given a radial cut profile $\lambda$, the simplicial fuzzy vector based on $\mathcal T$ and $\lambda$ has cuts
\begin{equation}\label{eq:simplicial}
 C_\alpha=\bmu+\lambda(\alpha)(\mathcal T-\bmu).
\end{equation}
Here $\lambda(\alpha)$ dilates the simplex about $\bmu$. Every vertex of $\mathcal T$ has Mahalanobis radius $\sqrt p$, so $C_\alpha$ lies in the Mahalanobis ellipsoid of radius $\sqrt p\,\lambda(\alpha)$ and meets its boundary at the vertices. The maps $R_{\bQ}$ commute with dilations about $\bmu$, so the completion \eqref{eq:completion} applies levelwise and turns $C_\alpha$ into exactly this ellipsoid. Symmetrizing the simplicial fuzzy vector thus yields $\Ell(\bmu,\bSigma,\rho)$ with
\begin{equation}\label{eq:lambdarho}
 \rho=\sqrt p\,\lambda.
\end{equation}

Its vertex selection is $\bX=\bmu+\lambda(U)(\bzeta_{N}-\bmu)$, with $N$ uniform and independent of $U$. Since $t(\bx)\le\sqrt p\lambda(\alpha)$ for $\bx\in C_\alpha$, the trace argument gives the sharp criterion
\[
 \text{simplicial carrying}\quad\Longleftrightarrow\quad\int_0^1\lambda(\alpha)^2\dd\alpha\ge1.
\]
The bound proves necessity, while scaling the vertex selection down when the integral exceeds one proves sufficiency. By \eqref{eq:lambdarho}, this is the ellipsoidal condition $\int_0^1\rho^2\ge p$ of Theorem~\ref{thm:carrying}, and simplicial normalization $\int_0^1\lambda^2=1$ corresponds to ellipsoidal normalization $\int_0^1\rho^2=p$.

The simplicial and the ellipsoidal carriers thus differ in their symmetry. The principal axes of an ellipsoidal carrier are those of $\bSigma$, so its orientation is fixed by the moment pair, and its cuts are invariant under every Mahalanobis rotation. The simplex of a sparsest realization, by contrast, has to point in a direction that the moment pair does not single out. For $p\ge2$, affine consistency rules out any such choice, although it leaves the radial profile free. Completing the simplex under all Mahalanobis rotations removes this additional orientation and produces the ellipsoid, which is the smallest crisp Mahalanobis ellipsoid that carries the moments.

Inclusion-minimality among all convex carriers does not single out the ellipsoid. For $p\ge2$ and $\int_0^1\lambda^2=1$, every cut of the simplicial fuzzy vector with $\lambda(\alpha)>0$ is a proper subset of the corresponding cut of $\Ell(\bmu,\bSigma,\sqrt p\,\lambda)$, although both fuzzy vectors carry the moments. Normalized ellipsoidal carriers are therefore minimal only within the fixed-shape Mahalanobis class, and it is affine consistency that singles out the ellipsoid. Similarly, every point of $\mathcal T$ has Mahalanobis radius at most $\sqrt p$, with equality only at the vertices. By \eqref{eq:trace}, every realizing law on a compact convex subset of $\mathcal T$ is therefore concentrated on the vertices, so that no proper compact convex subset of $\mathcal T$ carries the moments. Each regular simplex is thus a minimal crisp carrier. Distinct rotated copies of $\mathcal T$ are not nested, so for $p\ge2$ there is no smallest crisp carrier.

\begin{figure}
\centering
\begin{tikzpicture}[x=2.15cm,y=2.15cm,>=Stealth,font=\small]
 \draw[->] (-1.72,0)--(1.82,0) node[right] {$x_1$};
 \draw[->] (0,-1.67)--(0,1.74) node[above] {$x_2$};
 \foreach \x in {-1,1} {\draw (\x,0.025)--(\x,-0.025) node[below] {$\x$};}
 \foreach \y in {-1,1} {\draw (0.025,\y)--(-0.025,\y) node[left] {$\y$};}
 \draw[very thick,domain=0:360,samples=181]
   plot ({(3*cos(\x)+sin(\x))/sqrt(5)},{(cos(\x)+3*sin(\x))/sqrt(5)});
 \coordinate (A) at ({3/sqrt(5)},{1/sqrt(5)});
 \coordinate (B) at ({(-3+sqrt(3))/sqrt(20)},{(-1+3*sqrt(3))/sqrt(20)});
 \coordinate (C) at ({(-3-sqrt(3))/sqrt(20)},{(-1-3*sqrt(3))/sqrt(20)});
 \draw[thick] (A)--(B)--(C)--cycle;
 \foreach \pt in {A,B,C} {\fill (\pt) circle[radius=1.8pt];}
 \node[anchor=west,fill=white,inner sep=1pt] at (A) {$\bzeta_1$};
 \node[anchor=south east,fill=white,inner sep=1pt] at (B) {$\bzeta_2$};
 \node[anchor=north east,fill=white,inner sep=1pt] at (C) {$\bzeta_3$};
 \draw[dashed,thick]
   ({(3*cos(30)+sin(30))/sqrt(5)},{(cos(30)+3*sin(30))/sqrt(5)})--
   ({(3*cos(150)+sin(150))/sqrt(5)},{(cos(150)+3*sin(150))/sqrt(5)})--
   ({(3*cos(270)+sin(270))/sqrt(5)},{(cos(270)+3*sin(270))/sqrt(5)})--cycle;
\end{tikzpicture}
\caption{Sparsest realizations for the moments in Example~\ref{ex:bivariate}. The solid triangle joins three equiprobable points of Mahalanobis radius $\sqrt2$. The dashed triangle is another realization, rotated by $30^\circ$ in Mahalanobis coordinates. The union over all such rotations, not merely the two shown, is the outer ellipse $\{t\le\sqrt2\}$, the smallest carrying crisp Mahalanobis ellipsoid.}
\label{fig:simplex}
\end{figure}
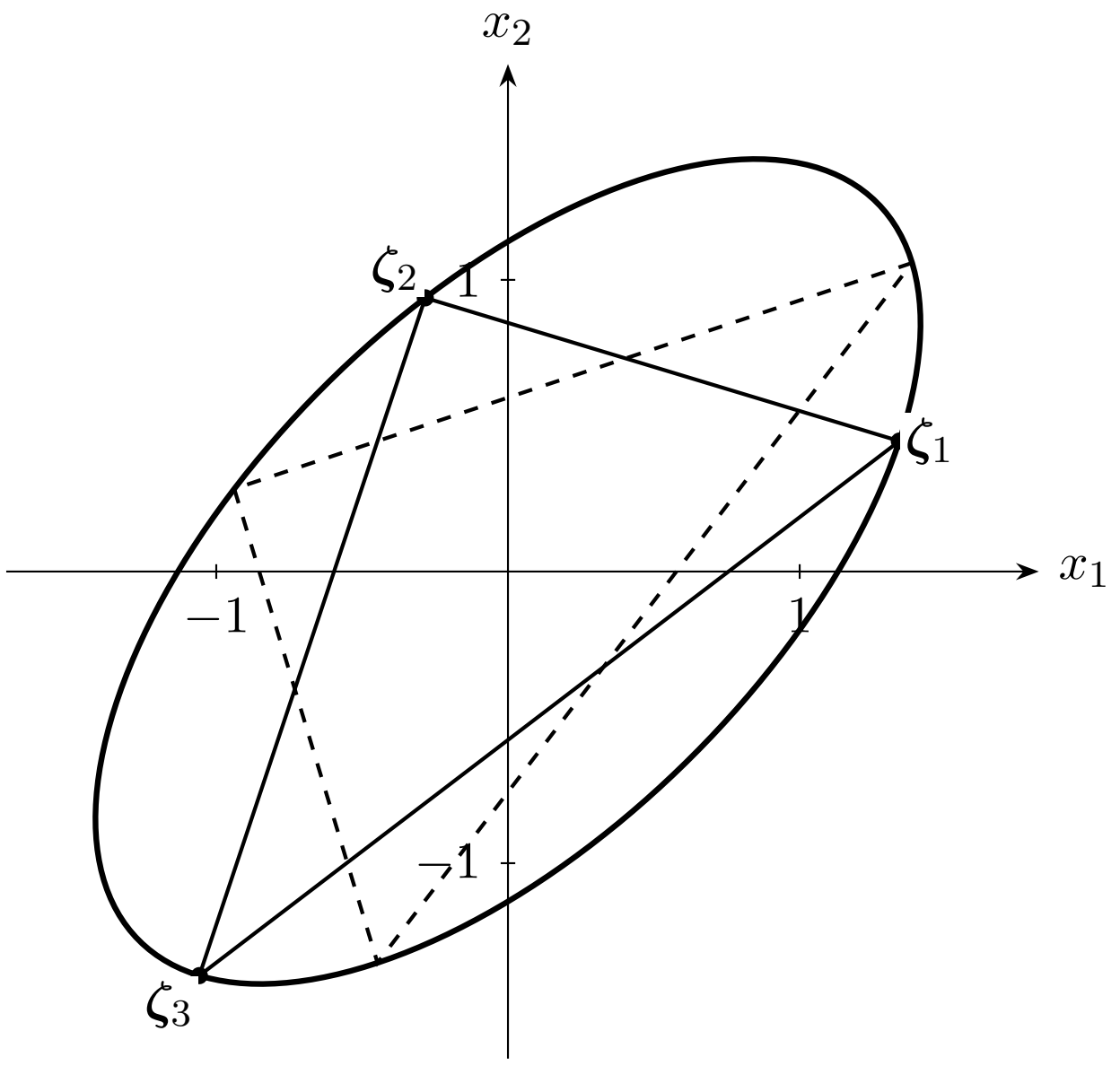

For illustration of the concepts of this section consider the following example.
\begin{example}[The same moments, a three-point realization]\label{ex:simplex}
Continue Example~\ref{ex:bivariate}, with $\bB=\bSigma^{1/2}$. Set
\[
 \bz_1=\begin{pmatrix}\sqrt2\\0\end{pmatrix},\quad
 \bz_2=\begin{pmatrix}-1/\sqrt2\\\sqrt{3/2}\end{pmatrix},\quad
 \bz_3=\begin{pmatrix}-1/\sqrt2\\-\sqrt{3/2}\end{pmatrix},\qquad
 \bzeta_j=\bB\bz_j.
\]
Directly, $\frac13\sum_j\bz_j=0$ and $\frac13\sum_j\bz_j\bz_j^\top=\bI_2$, so the three equally weighted $\bzeta_j$ realize $(0,\bSigma)$. They all satisfy $t(\bzeta_j)=\sqrt2$, and their convex hull $\mathcal T$ is a crisp carrier.
Its Mahalanobis-rotation orbit fills the ellipse $\{t\le\sqrt2\}$ in Figure~\ref{fig:simplex}. The construction selects one triangle only to display a realization. Another rotated triangle is equally compatible with the moments.

For the two profiles of Example~\ref{ex:bivariate}, the corresponding normalized simplicial profiles $\lambda=\rho/\sqrt p$ are $\lambda_{\mathrm{tri}}(\alpha)=\sqrt3(1-\alpha)$ and $\lambda_{\mathrm{par}}(\alpha)=\sqrt{2(1-\alpha)}$.

Their symmetrized cuts recover the two ellipsoidal fuzzy vectors. The example therefore separates three choices. The affine consistency conditions~\ref{prop:coordinateneut} and \ref{prop:ambientneut} of Definition~\ref{def:affine} force the ellipsoidal geometry, carrying restricts the radial scale, and a discrete realization requires an orientation that the moments do not specify.
\end{example}

\section{Conclusion}\label{sec:discussion}

We have studied how to assign a fuzzy vector to a prescribed mean vector and covariance matrix without assuming an underlying probability law. Starting from requirements on the representation rather than from a membership function, we showed that, for fuzzy vectors with compact convex cuts, coordinate and ambient neutrality hold if and only if the cuts are homothetic Mahalanobis ellipsoids, with one radial profile for each positive covariance rank. The moments thus fix the geometry, and the radial profile remains the modeling choice.

We introduced \emph{carrying}, the existence of a random-set selection with exactly the prescribed moments. For ellipsoidal fuzzy vectors it is equivalent to a sharp integral bound on the squared radial profile, and equality characterizes the non-unique inclusion-minimal carriers within the fixed-shape Mahalanobis class. Compatibility also with rank-reducing maps would force a common profile across ranks. Within the bounded-support class considered here, this is incompatible with carrying moment pairs of every positive rank, and with normalization at two distinct positive ranks. Sparsest equally weighted realizations of the moments are regular simplices in Mahalanobis coordinates, and their rotational completion recovers the unique minimum-volume crisp ellipsoid in the state space that carries the moments. The ellipsoidal representation finally encodes interaction through the covariance quadratic form.

\paragraph{Scope and further research}
 Neither ellipsoids nor homothetic cuts are properties of arbitrary carriers. For example, nested ellipsoidal cuts may have level-dependent centers and shape matrices, $C_\alpha=\bc(\alpha)+\bL(\alpha)\mathbb B_d$. Such fuzzy vectors are legitimate when they satisfy the cut conditions, but this additional level-dependent structure cannot arise from an affine-consistent construction based only on $(\bmu,\bSigma)$. It must be elicited, modeled or inferred from additional information, or obtained by relaxing the neutrality requirements.

The present theory treats the moments as given. It does not propagate estimation error in $\bSigma$, identify skewness or nonlinear dependence, or recover tail behavior from covariance alone. The ellipsoidal joint is a model satisfying the axioms, not a reconstruction of an unknown joint probability law. Likewise, exact treatment of singular matrices is not a claim of continuity when the covariance rank changes, because rank-indexed profiles may change at such a transition.

There are three distinct directions for further work. One is to identify principles weaker than compatibility with all affine maps, which Proposition~\ref{prop:rankreduction-carrying} rules out, that relate profiles across intrinsic ranks. Another is to recover the moment pairs admitted by a specified fuzzy vector, rather than construct a vector from prescribed moments. This selection-based inverse question differs from defining a mean and dispersion directly from the membership function, an approach explored for fuzzy numbers in \cite{Sch20}. A third direction is to admit higher-order information, such as third moments. With such information, the admissible symmetries are those preserving the additional tensor, and the compatible cut geometries would have to be analyzed accordingly, possibly with the help of an additional selection principle. The truncated moment problem \cite{CurtoFialkow2005,Lasserre2010} and work on skewness \cite{Sch20,JKPT} offer starting points.

\section*{Acknowledgments}
Paul Schneider gratefully acknowledges financial support from the Swiss National Science Foundation (SNSF) through the project ``Large-scale kernel methods in financial economics'' (grant no.~215528).

\bibliographystyle{elsarticle-num}
\bibliography{references}
\end{document}